\documentstyle[12pt]{article}

  \newcommand{\const}{\rm const}
  \newcommand{\Var}{\rm Var}

  \newcommand{\vraisup}{\rm vraisup}
  \newcommand{\Dom}{\rm  Dom}

  \newcommand{\meas}{\rm meas}
  \newcommand{\Sub}{\rm  Sub}

\begin{document}

   \begin{center}

{\bf  Low bounds for distribution of sums  of independent centered } \\

\vspace{4mm}

{\bf random variables belonging to  Grand Lebesgue Spaces. }\\

  \vspace{4mm}

 {\bf M.R.Formica, \ E.Ostrovsky, \ L.Sirota. }

\vspace{3mm}

\end{center}

 \ Universit\`{a} degli Studi di Napoli Parthenope, via Generale Parisi 13, Palazzo Pacanowsky, 80132,
Napoli, Italy. \\

e-mail: mara.formica@uniparthenope.it \\

\vspace{3mm}

 \ Israel,  Bar - Ilan University, department  of Mathematic and Statistics, 59200. \\
e-mails: \  eugostrovsky@list.ru \\
sirota3@bezeqint.net \\

\vspace{4mm}

 \begin{center}

 {\bf Abstract} \\

 \vspace{3mm}

 \ We  deduce in this  short report the non - asymptotic  {\it lower bounds} for exponential tail of distribution for sums of independent
 centered random variables.\\

 \end{center}

 \vspace{5mm}

\begin{center}

 {\it Key words and phrases.} \\

\end{center}

\vspace{4mm}

 \  Probability space, centered random variable (r.v.),   Lebesgue - Riesz and Grand Lebesgue Spaces (GLS) and norms, natural function,  independence,
 tail of distribution, generating function,  rearrangement invariant Banach functional spaces, anti - norm, anti - triangle inequality,
  upper and lower estimates. \par

\vspace{5mm}

\section{Definitions. Notations. Statement of problem.}

\vspace{4mm}

 \ Let $ \  (\Omega, B, {\bf P}) \ $ be certain probability space with expectation $ \ {\bf E} $  and dispersion  $ \ \Var;  \ X, Y \ $ be independent centered (mean zero) random variables (r.v.)
 and $ \ q = \const \in [1, \infty]. \ $ The ordinary  Lebesgue - Riesz, or $ \ L(q) \ $ norm of the  arbitrary r.v. $ \ Z \ $  will be denoted by $ \ |Z|_q: \ $

$$
|Z|_q := \left[ \ {\bf E} |Z|^q \ \right]^{1/q}, \ 1 \le q < \infty,
$$
and

$$
 |Z|_{\infty} := \vraisup_{\omega \in \Omega} |Z(\omega)|.
$$

 \  Assaf Naor and Krzysziof Oleszkiewicz in a recent article \cite{Naor} proved in particular the following inequality for the r.v.- s.
 belonging to some Lebesgue - Riesz space

\begin{equation} \label{Naor in}
|X + Y|_q  \ge  \left[ \ |X|^q_q + |Y|^q_q  \  \right]^{1/q}, \ q \in [2,\infty],
\end{equation}
in our notations. More generally,  let $ \ \{ X_i \}, \ i = 2,3, \ldots,n  \ $ be a family of (common) independent centered r.v. - s; then by induction

\begin{equation} \label{n case q}
 \left| \ \sum_{i=1}^n X_i \   \right|_q  \ge \left[ \ \sum_{i=1}^n  |X_i|^q_q  \right]^{1/q}, \ q \in [2,\infty].
\end{equation}
 \ If in addition the r.v. $ \ X_i \ $ are identical distributed,

\begin{equation} \label{n  norm case q}
 \left| \ n^{-1/2} \sum_{i=1}^n X_i \   \right|_q  \ge  n^{1/q - 1/2} \  |X_1|_q, \ q \in [2,\infty].
\end{equation}

 \vspace{3mm} 

\ The estimation (\ref{n  norm case q}) may be named as {\it power level}.

 \ Note that this result is weak if $ \ q > 2 \ $ as $ \ n \to \infty;\ $  later we will improve it.\par

\vspace{5mm}

 \ {\bf Our  purpose in this short article is to extend the last inequality into the r.v. belonging the so - called
 Grand Lebesgue Spaces (GLS).  } \par

 \vspace{5mm}
 
 \ We obtain as a consequence an exact non - uniform {\it lower} exponential estimations for tail of distribution for the sums of
independent centered r.v. \par

\vspace{4mm}

\begin{center}

\ {\sc Brief note about  Grand Lebesgue Spaces (GLS). } \\

\end{center}

{\it A classical approach.} \par

 \vspace{3mm}

 \ Let $\lambda_0 \in (0,\infty]$ and let $ \phi = \phi(\lambda)$ be an
even strong convex function in $(-\lambda_0, \lambda_0)$ which takes
positive values, twice continuously differentiable; briefly $\phi =
\phi(\lambda)$ is a Young-Orlicz function, such that

\begin{equation}\label{Young-Orlicz function}
\phi(0) = 0, \ \ \phi'(0) = 0, \ \  \phi^{''}(0) \in (0,\infty).
\end{equation}

 \ We denote the set of all these Young-Orlicz function as  $\Phi: \ \Phi = \{ \phi(\cdot)  \}. $ \\

{\bf Definition 1.1.}\par

 \ Let $\phi\in \Phi$. We say that the centered random variable $\xi$
belongs to the space $B(\phi)$  if there exists a constant $\tau \geq 0$ such that

\begin{equation}\label{spaceB}
\forall \lambda \in (-\lambda_0, \lambda_0) \ \Rightarrow {\bf E}
\exp(\pm \lambda \ \xi) \le \exp(\phi(\lambda \ \tau)).
\end{equation}

 \ The minimal non-negative value $\tau$ satisfying ( \ref{spaceB}) for
any $\lambda \in (-\lambda_0, \ \lambda_0)$ is named $B(\phi)$-norm
of the variable $\xi$  and we write

\begin{equation}\label{Bnorm}
||\xi||_{B(\phi)} \stackrel{def}{=}\inf \{\tau \ge 0 \ : \ \forall
\lambda \in (-\lambda_0, \lambda_0) \ \Rightarrow {\bf E} \exp(\pm
\lambda \ \xi) \le \exp(\phi(\lambda \ \tau)) \} .
\end{equation}

\vspace{3mm}

 \ For instance if $\phi(\lambda)=\phi_2(\lambda) := 0.5 \ \lambda^2, \
\lambda \in \mathbf{R}$, the r.v. $\xi$ is \emph{subgaussian} and in
this case we denote the space $B(\phi_2)$ with $\Sub$. Namely we
write $\xi \in \Sub$ and
$$
||\xi||_{\Sub} \stackrel{def}{=} ||\xi||_{B(\phi_2)}.
$$
 \  It is known, see  \cite{KozOs}, \cite{Buldygin-Mushtary-Ostrovsky-Pushalsky} that if the r.v. $\xi_i$ are
independent and subgaussian, then

\begin{equation} \label{Sums estim}
||\sum_{i=1}^n \xi_i||_{\Sub} \le \sqrt{\sum_{i=1}^n ||\xi_i||^2_{\Sub}}.
\end{equation}

 \ At the same inequality holds true in the more general case in the $ \ B(\phi) \ $ norm, when the function $ \ \lambda \to \phi(\sqrt{\lambda})  \ $ is convex,
 see \cite{KozOs}. \par
 \ As a slight corollary: in this case and if in addition the r.v. - s $ \ \{\xi_i \} \ $ are i., i.d., then

\begin{equation} \label{upp for sums}
 \sup_{n = 1,2,\ldots}|| n^{-1/2} \sum_{i=1}^n \xi_i||B(\phi) =  ||\xi_1||B(\phi).
\end{equation}

\vspace{4mm}


\vspace{3mm}

It is proved in particular that $B(\phi), \ \phi  \in \Phi$, equipped with the norm
(\ref{Bnorm}) and under the ordinary algebraic operations, are
Banach rearrangement invariant  functional spaces, which are
equivalent the so-called Grand Lebesgue spaces as well as to Orlicz
exponential spaces. These spaces are very convenient for the
investigation of the r.v. having an exponential decreasing tail of
distribution; for instance, for investigation of the limit theorem,
the exponential bounds of distribution for sums of random variables,
non-asymptotical properties, problem of continuous and weak
compactness of random fields, study of Central Limit Theorem in the
Banach space, etc. \par

\vspace{3mm}

 \ Let $ \ g: R \to R \ $ be numerical valued measurable function, which can perhaps take the infinite value.
 Denote by $ \Dom[g] \ $ the domain of its finiteness:

\begin{equation}\label{Domain}
\Dom[g] := \{y, \ g(y) \in (-\infty, \ + \infty) \ \}.
\end{equation}

 \   Recall the definition $ \ g^*(u) \ $  of  the Young-Fenchel or
Legendre transform for the function $ \ g: R \to R  \ $:

\begin{equation}\label{ Definition of the Young-Fenchel transform}
g^*(u) \stackrel{def}{=} \sup_{y \in \Dom[g]} (y u - g(y)),
\end{equation}
 but we will use further the value $ \ u \ $  to be only non - negative.\par

 \ In particular, we  denote by $\nu(\cdot)$ the Young-Fenchel or
Legendre transform for the function $\phi\in \Phi$:

\begin{equation}\label{Young-Fenchel transform}
\nu(x) = \nu[\phi](x)  \stackrel{def}{=} \sup_{\lambda: |\lambda|
\le \lambda_0} (\lambda x - \phi(\lambda)) = \phi^*(x).
\end{equation}

 \ It is important to note that if the non-zero r.v. $\xi$ belongs to
the space $B(\phi)$ then

\begin{equation}\label{conditionB}
{\bf P}(\xi > x) \le \exp \left(- \nu(x/||\xi||_{B(\phi)}\right).
\end{equation}
 \ The inverse conclusion is also true up to a multiplicative constant
under  suitable conditions.\par

\vspace{4mm}

 \ Furthermore, assume that the {\it centered} r.v. $\xi$  has in some
non-trivial neighborhood of the origin finite \emph{moment generating
function} and define

\begin{equation}\label{momentfunction}
\phi_{\xi}(\lambda) \stackrel{def}{=} \max_{\alpha = \pm 1} \ \ln {\bf E} \exp (\ \alpha \lambda \ \xi \ ) < \infty, \ \lambda \in ( - \ \lambda_0, \lambda_0)
\end{equation}
for some $\lambda_0 = \const \in (0, \ \infty]$. Obviously, the last
condition (\ref{conditionB}) is quite equivalent to the well known
Cramer's one. \par

 \ We agree that $\phi_{\xi}(\lambda) := \infty$ for all the values $\lambda$ for which
\begin{equation}\label{mean infinity}
 {\bf E} \exp ( \ |\lambda| \ \xi) = \infty.
\end{equation}
The function $\phi_{\xi}(\lambda)$ introduced in
(\ref{momentfunction}) is named {\it natural} function for the r.v.
$ \xi$; herewith  $\xi \in B(\phi_{\xi}) $ and moreover we assume

$$
||\xi||_{B(\phi_{\xi})} = 1.
$$

\vspace{5mm}

 \ We recall here for reader convenience some known definitions and  facts about  Grand Lebesgue Spaces (GLS) using in this article.

   \ Let $ \ \psi = \psi(p), \ p \in [1,b) \ $    where $ \ b = \const, \ 1 \le   b \le \infty \ $ be positive measurable numerical valued
    function, not necessary to be finite in every point, such that $ \ \inf_{p \in [1,b)} \psi(p) > 0. \ $   For instance

  $$
    \psi_m(p) := p^{1/m}, \ m = \const > 0, \ p \in [1,\infty)
  $$
  or

$$
   \psi^{(b; \beta)}(p) :=  (b-p)^{-\beta}, \ p \in [1,b), \ b = \const, \  1 \le b < \infty; \ \beta = \const \ge 0.
$$

\vspace{4mm}

{\bf Definition 1.2.} \par

\vspace{3mm}

 \ By definition, the (Banach) Grand Lebesgue Space (GLS)    $  \ G \psi  = G\psi(b),  $
    consists on all the real (or complex) numerical valued random variable (measurable functions)
   $   \  f: \Omega \to R \ $  defined on whole our  space $ \ \Omega \ $ and having a finite norm

 \begin{equation} \label{norm psi}
    || \ f \ || = ||f||G\psi \stackrel{def}{=} \sup_{p \in [1,b)} \left[ \frac{|f|_p}{\psi(p)} \right].
 \end{equation}

 \vspace{4mm}

 \ The function $ \  \psi = \psi(p) \  $ is named as  the {\it  generating function } for this space. \par

  \ If for instance

$$
  \psi(p) = \psi^{(r)}(p) = 1, \ p = r;  \  \psi^{(r)}(p) = +\infty,   \ p \ne r,
$$
 where $ \ r = \const \in [1,\infty),  \ C/\infty := 0, \ C \in R, \ $ (an extremal case), then the correspondent
 $ \  G\psi^{(r)}(p)  \  $ space coincides  with the classical Lebesgue - Riesz space $ \ L_r = L_r(\Omega, {\bf P}). \ $ \par

\vspace{4mm}

 \ These spaces are investigated in many works, e.g. in
 \cite{Fiorenza1},   \cite{Fiorenza3}, \cite{Fiorenza4},   \cite{Iwaniec1}, \cite{Iwaniec2}, \cite{KozOs},
\cite{LiflOstSir},   \cite{Ostrovsky1}  - \cite{Ostrovsky5} etc. They are applied for example in the theory of Partial Differential Equations
\cite{Fiorenza3}, \cite{Fiorenza4}, in the theory of Probability  \cite{Ermakov},\cite{Ostrovsky3}  - \cite{Ostrovsky5}, in Statistics \cite{Ostrovsky1}, chapter 5,
theory of random fields  \cite{KozOs}, \cite{Ostrovsky4}, in the Functional Analysis \cite{Ostrovsky1}, \cite{Ostrovsky2}, \cite{Ostrovsky4} and so one. \par

 \  These spaces are rearrangement invariant (r.i.) Banach functional spaces; its fundamental function  is considered in  \cite{Ostrovsky4}. They
  not coincides  in general case with the classical  spaces: Orlicz, Lorentz, Marcinkiewicz  etc., see \cite{LiflOstSir} \cite{Ostrovsky2}.

 \  The belonging of  some  r.v. $ \ f:  \Omega \to R \ $ to some $ \ G\psi \ $ space   is closely related with its tail behavior

  $$
  T_f(t) = \meas \left\{x; \ x  \in R^d, \ |f(x)| > t  \right\}
  $$
   as $ \ t \to \infty,  \ $  see  \cite{KozOs}, \cite{KozOsSir2017}.    \par

  \ Let a family of the functions $ \ \{ f_w   \} =  \{ f_w(\omega) \}, \ x \in R^d, \ w \in W, \  $   where $ \ W = \{w\} \ $ is arbitrary set, be such that

 \begin{equation} \label{cond psi}
 \exists b \in [1, \infty] \ \Rightarrow \psi[W](p) := \sup_{ p \in (a,b) } |f|_p < \infty.
 \end{equation}

 \ The function $ \ \psi[W](p) \ $ is named as a {\it natural function}   for the family $ \  \{ f_w   \}, \ w \in W. $  It may be considered as a generating function for certain
Grand Lebesgue Space $ \ G\psi[W]. \ $ Obviously,

$$
\sup_{w \in W} ||f_w||G\psi[W] = 1.
$$
 \ Notice that the family $ \ \{f_w\}  \ $ may consists on the single function $ \ f_w = f, \ $ if course it satisfied the condition (\ref{cond psi}); we will write then

 $$
\psi[f](p) := |f|_p, \  1 \le p  < b,
 $$
one can take $ \ 1 \le p \le b, \ $ if $ \ B < \infty \ $ and $ \ |f|_b < \infty. \ $  \par

 \vspace{5mm}

\section{Main result: lower estimate. Anti - norms.}

\vspace{4mm}

 \ The theory of GLS allows in particular to deduce the {\it upper } bound for distribution of sums of random variables,
independent, centered or not. The norm in these spaces is defined by means of the  operation $ \ \sup, \ $ see (\ref{norm psi}). \par
 \ It is reasonable to assume that for an obtaining of the {\it lower } bounds  for these sums we must apply for definition of some
 functionals of a type "norm"  use the operator $  \  \inf. \ $  In detail:\par

 \vspace{4mm}

\ {\bf Definition 2.1.}  Let $ \ \psi  = \psi(p), \ p \in [1.b) \ $ be certain generating function: \ $ \psi(\cdot) \in \Psi_b. \ $ The following
functional is named as the  $ \ AG\psi - $ {\it anti - norm }  \\
 $ \ V(X) = V(X)AG\psi = V(X)_{\psi}\ $ of the r.v. $ \ X: \ $

 \vspace{4mm}

\begin{equation} \label{anti norm psi}
 V(X) = V(X)AG\psi  \stackrel{def}{=} \inf_{p \in [1,b)} \left[ \frac{|X|_p}{\psi(p)} \right],
\end{equation}
in contradiction with the classical definition of the GLS norms, see (\ref{norm psi}). \par

 \ The (linear) space of all the random variables having non trivial  $ \ AG\psi - $ norms forms by definition the Anti - Grand Lebesgue space $ \   AG\psi. \ $ \par

\vspace{4mm}

 \ The following  properties of introduced anti - norm are evident: $  \ V(X) \ge 0; \ $ and if in addition the generating function  $ \ \psi(\cdot) \ $ is bounded:
$ \ \sup_p \psi(p) < \infty, \ $ then

$$
V(X) \ge 0, \ V(X) =  0 \Leftrightarrow X = 0;
$$

$$
\forall C \in R \ \Rightarrow  V(CX) = |C| \ V(X);
$$

$$
V(X + Y) \ge V(X) + V(Y) \ -
$$
anti - triangle inequality. \par

\vspace{3mm}

 \ {\bf Remark 2.1.}  If the generating function $ \ \psi(\cdot) \ $ coincides with the natural function of some
 r.v. $ \ X, \ \psi(p) = \psi[X](p) =  |X|_p, \ 1 \le p < b, \ $  then obviously the ordinary and anti- GLS norms  of r.v. $ \ X \ $ coincides:

$$
||X||G\psi = ||X|| AG\psi = V(X)  = 1.
$$

\vspace{3mm}

 \ Let us now investigate the strengthening of the anti - triangle inequality for independent centering r.v.  $ \ X \ $ and $ \ Y. \ $ Define for this purpose
the following functions

$$
\theta(p,q) := \inf_{a,b > 0} \left[ \ \frac{(a^q + b^q)^{1/q}}{(a^p + b^p)^{1/p}} \ \right] =
\inf_{z > 0} \left[ \ \frac{(z^q + 1)^{1/q}}{(z^p + 1)^{1/p}} \ \right], \ p,q \ge 1;
$$

$$
\kappa(p) = \kappa_b(p)  := \min_{q \in [1,b)} \theta(p,q);
$$

then

\begin{equation} \label{low a b}
 \forall p \ge 1 \ \Rightarrow   (a^q + b^q)^{1/q} \ge \theta(p,q) \ (a^p + b^p)^{1/p},
\end{equation}

\begin{equation} \label{theta}
\theta(p,q) = \min \left(\ 1, 2^{1/q - 1/p}    \ \right).
\end{equation}

\begin{equation} \label{kappa b}
 \kappa_b(p) = \min \left( \ 1, 2^{1/b - 1/p}  \  \right),
\end{equation}
so that

\begin{equation} \label{kappab  p le b}
 \kappa_b(p) =  2^{1/b - 1/p},   \  1 \le p \le b,
\end{equation}
 and

\begin{equation} \label{kappab p ge b}
 \kappa_b(p) = 1,  \  p > b.
\end{equation}

\vspace{5mm}

 \ Let now the centered independent random variables $ \ X, \ Y \ $ belongs to some
Anti - Grand Lebesgue space $ \   AG\psi, \ \exists \psi \in \Psi(b): \ $

$$
|X|_q \ge V(X) \psi(q), \ |Y|_q \ge V(Y) \psi(q), \ 1 \le q < b.
$$
 \ We apply the  Naor and Oleszkiewicz inequality (\ref{Naor in})  for the values $ \  q \in [1, b) \ $

\begin{equation} \label{Naor  psi}
|X + Y|_q  \ge  \left[ \ |X|^q_q + |Y|^q_q  \  \right]^{1/q} \ge \psi(q) \ \left[ \ |V(X)|^q + |V(Y)|^q \  \right]^{1/q};
\end{equation}

$$
\frac{|X + Y|_q }{\psi(q)}  \ge  \left[ \ |V(X)|^q + |V(Y)|^q \  \right]^{1/q}.
$$
 \ Let now $ \ p = \const \ge 1; \ $  we obtain using (\ref{low a b}) and  (\ref{kappa b})

\begin{equation} \label{V X Y p}
V(X + Y) \ge \kappa_b(p) \ \left (V^p(X) + V^p(Y) \ \right)^{1/p}, \ p \ge 1.
\end{equation}

 \ Highlight a particularly very important case $ \ p = 2: \ $

\begin{equation} \label{V X Y 2}
V(X + Y) \ge \kappa_b(2) \ \left (V^2(X) + V^2(Y) \ \right)^{1/2}.
\end{equation}

  \ More detail:

\begin{equation} \label{V X Y  detail}
V(X + Y) \ge \min  \left( 1, 2^{1/b - 1/p}   \right)  \ \left (V^p(X) + V^p(Y) \ \right)^{1/p}.
\end{equation}

 \vspace{3mm}

 \ To summarize: \par

 \vspace{5mm}

{\bf Theorem 2.1.} Let $ \ \{X_i\}, \ i = 1,2,\ldots,n \ $ be a sequence  of centered independent random variables  belonging to some
Anti - Grand Lebesgue space  $ \ AG\psi, \ \psi \in \Psi(b), \ 1 < b \le \infty. \ $  Our proposition:

\begin{equation} \label{theorem}
V \left( \ \sum_{i=1}^n X_i \  \right) \ge \min \left(1, \ 2^{1/b - 1/p}   \right) \left[ \ \sum_{i=1}^n V^p(X_i) \  \right]^{1/p}, p \in [1,\infty].
\end{equation}

\ In particular:

\begin{equation} \label{coroll p 2}
V \left( \ \sum_{i=1}^n X_i \  \right) \ge \min \left(1, \ 2^{1/b - 1/2}   \right) \left[ \ \sum_{i=1}^n V^2(X_i) \  \right]^{1/2}.
\end{equation}

 \ If in addition $ \ b = \infty, \ $ then

\begin{equation} \label{coroll b infty}
V \left( \ \sum_{i=1}^n X_i \  \right) \ge  2^{- 1/2}   \ \left[ \ \sum_{i=1}^n V^2(X_i) \  \right]^{1/2}.
\end{equation}

\vspace{5mm}

 \section{ Examples.}

\vspace{4mm}

 {\bf A.} Let us consider the symmetrical distributed  subgaussian r.v $ \ X \ $ defined on some sufficiently rich probability space having the density

\begin{equation} \label{ex A}
f_X(x) = 0.5 \ |x| \ e^{-x^2/2}, \ x \in (-\infty, \infty).
\end{equation}

 \ We have for non - negative values $ \ p \ $

$$
{\bf E} |X|^p = 2^{p/2} \ \Gamma(p/2 + 1),
$$
therefore the natural function for this r.v. is following

\begin{equation} \label{psi examp}
\psi[X](p) = |X|_p = 2^{1/2} \ [\Gamma(p/2 + 1) ]^{1/p}.
\end{equation}

 \ Note that as $ \ p \in [1,\infty) \ $

$$
\psi_X(p) \asymp (p/e)^{1/2}.
$$

\vspace{4mm}

 {\bf B.} Let the centered r.v. $ \ X \ $  be {\it bilateral subgaussian:}

$$
C_1 p^{1/2} \le \psi[X](ç) \le C_2 p^{1/2}, \ \exists C_1, C_2 \in (0,\infty), \ C_1 \le C_2,  \ 1 \le p < \infty.
$$
 \ Let also $ \ X_i, \ i = 1,2,\ldots \ $ be independent copies of $ \ X. \ $ Define the classical normed sum

$$
S_n := n^{-1/2} \sum_{i=1}^n X_i.
$$

 \ We deduce by virtue of theorem 2.1 that for  $ u \ge 1 \ $

 $$
 \exists C_3, C_4 = \const \in (0,\infty), 0 < C_4 \le C_3 \ \Rightarrow
 $$

\begin{equation} \label{Subg est}
 \exp (-C_3 u^2) \le {\bf P} (S_n >u) \le \exp (-C_4 u^2)
\end{equation}

and  and the same estimate there holds  for left - hand side tail $ \ {\bf P} (S_n < -u). \ $ \\

\vspace{4mm}

{\bf C.} \  Let us consider a more general case of the sequence of centered independent  r.v. $ \ \{X_1, X_2, \ldots,X_n\} \ $ such that

$$
\exists m > 0, \ \exists C_5,C_6 \in (0,\infty), C_6 \le C_5, \ \forall u \ge 1 \ \Rightarrow
$$

$$
\exp(-C_5 u^m ) \le {\bf P}(|X_i| > u) \le \exp(-C_6 u^m),
$$
or equally

$$
C_7 \ p^{1/m} \le \inf_i \psi[X_i](p) \le \sup_i \psi[X_i](p) \le  C_8 \ p^{1/m}, \ p \in [1,\infty).
$$
 \ We propose

$$
\exists C_9,C_{10} \in (0,\infty), \ C_{10} \le C_9 \ \Rightarrow \ \exp \left(-C_9  u^{\min(m,2)}  \ \right) \le
$$

\begin{equation} \label{min m 2}
 {\bf P} \left( \ n^{-1/2} \left|\sum_{i=1}^n X_i \right| > u \ \right) \le
\exp \left(-C_{10} \  u^{\min(m,2)}  \ \right), \ u \ge 1.
\end{equation}

 \ Note that the upper estimate in (\ref{min m 2} ) is known, see \cite{KozOs}, \cite{Ostrovsky1}, chapter 2, section 2.1.\par

\vspace{5mm}

\vspace{0.5cm} \emph{Acknowledgement.} {\footnotesize The first
author has been partially supported by the Gruppo Nazionale per
l'Analisi Matematica, la Probabilit\`a e le loro Applicazioni
(GNAMPA) of the Istituto Nazionale di Alta Matematica (INdAM) and by
Universit\`a degli Studi di Napoli Parthenope through the project
\lq\lq sostegno alla Ricerca individuale\rq\rq (triennio 2015 - 2017)}.\par

\end{document}